\documentclass[12pt,reqno]{amsart}

\usepackage{color}
\usepackage{esint,amssymb}
\usepackage{graphicx}
\usepackage{tikz}

\usepackage{hyperref}

\usepackage{amsmath}
\numberwithin{equation}{section}
\usepackage{mathtools}
\mathtoolsset{showonlyrefs}

\newtheorem{theorem}{Theorem}
\newtheorem{proposition}[theorem]{Proposition}

\theoremstyle{definition}

\numberwithin{equation}{section}
\numberwithin{theorem}{section}

\newcommand{\e}{{\varepsilon}}

\newcommand{\R}{\mathbb{R}}
\renewcommand{\S}{\mathbb{S}}

\newcommand{\ep}{\varepsilon}

\newcommand{\average}{{\mathchoice {\kern1ex\vcenter{\hrule height.4pt
width 6pt depth0pt} \kern-9.7pt} {\kern1ex\vcenter{\hrule
height.4pt width 4.3pt depth0pt} \kern-8pt} {} {} }}
\newcommand{\ave}{\average\int}

\begin{document}

\title[Stable solutions to reaction-diffusion elliptic problems]
{Stable solutions to reaction-diffusion elliptic problems}

\author[X. Cabr\'e]{Xavier Cabr\'e}
\address{﻿X.C.\textsuperscript{1,2,3} 
\textsuperscript{1}ICREA, Pg.\ Lluis Companys 23, 08010 Barcelona, Spain \& 
\textsuperscript{2}Universitat Polit\`ecnica de Catalunya, Departament de Matem\`{a}tiques and IMTech, 
Diagonal 647, 08028 Barcelona, Spain \&
\textsuperscript{3}Centre de Recerca Matem\`atica, Edifici C, Campus Bellaterra, 08193 Bellaterra, Spain
}
\email{xavier.cabre@upc.edu}

\begin{abstract}

We are concerned with stable solutions to reaction-diffusion elliptic PDEs. We begin with regularity questions, first addressing the classical Laplacian. In joint work with Figalli, Ros-Oton, and Serra~\cite{CFRS}, we proved that stable solutions are smooth up to the optimal dimension 9, thereby solving an open problem posed by Brezis in the mid-1990s. We describe this result and also discuss related progress and open problems for the fractional Laplacian ---arising naturally in boundary reaction problems---, the $p$-Laplacian, and minimal surfaces.

We then turn to existence questions, starting with the Casten-Holland and Matano theorem for interior reactions, which states that no nonconstant stable solution exists in convex domains  under zero Neumann boundary conditions. We present a recent result with C\'onsul and Kurzke~\cite{CCK} establishing that the analogous statement fails for boundary reactions. This requires the development of a new Ginzburg-Landau theory for real-valued functions and the analysis of the half-Laplacian on the real line, for which we present new results and open problems.

\end{abstract}


\thanks{The author was supported by grants PID2021-123903NB-I00 and RED2022-134784-T funded by MCIN/AEI/10.13039/501100011033 and by ERDF ``A way of making Europe'', and by the Catalan grant 2021-SGR-00087, as well as by the Spanish State Research Agency through the Severo Ochoa and Mar\'{\i}a de Maeztu Program for Centers and Units of Excellence in R\&D (CEX2020-001084-M)}

\maketitle

\section{Hilbert's 19th problem on regularity}

In many physical and geometric problems, the states that can be observed are those that minimize a certain functional. When the possible states are described by functions $u$ of one or several real variables, the functional is a real-valued map $\mathcal E$ acting on such functions. In classical mechanics, $\mathcal E$ is called the {\it action} and is given by the integral of a Lagrangian. A simple example is the motion of a particle under gravity: its position is represented by $u=u(x)$ with $x=t\in \R$ denoting time, and the action is the integral of kinetic minus potential energies. In geometry, two fundamental examples are geodesics ---curves in a Riemannian manifold that are critical points of the length functional--- and minimal surfaces, which are hypersurfaces in Euclidean space that are critical points of the area functional.  

Hilbert’s 19th problem asks whether minimizers of elliptic functionals are necessarily analytic. For functionals of the form  
\[
\mathcal E(u)=\int_{\Omega} L(\nabla u(x))\,dx,
\]  
with $L:\R^n\to\R$ a uniformly convex function and $u:\Omega\subset \R^n\to\R$, this was solved in the late 1950s through the independent work of Ennio De Giorgi and John Forbes Nash, Jr. They proved a H\"older continuity result that answers in the affirmative Hilbert's question in case $L$ is analytic. Our work \cite{CFRS}, described below, revisits the same question for Lagrangians of the form  
\[
L(\nabla u, u)=\tfrac{1}{2}|\nabla u|^2 - F(u),
\]  
which depend explicitly on the variable $u$.  

The {\it principle of least action} in mechanics states that observable states are not merely critical points of the action, but absolute minimizers. In the convex case ---when $L=L(\nabla u)=L(p)$ is a convex function of $p\in\R^n$--- the functional $\mathcal E$ is itself convex, and hence any critical point (if it exists) is automatically an absolute minimizer. In this setting, the principle of least action is perfectly valid.  

However, in some real-world situations the principle fails since the associated functional is not convex. We may observe states that are only {\it local minimizers} ---minimizers among small perturbations--- even in cases where an absolute minimizer also exists. Such states are, in particular, {\it stable states}, meaning that the second variation of the functional at these points is nonnegative definite.  Below, we will give two concrete real-world examples of this fact.

\section{Stable minimal surfaces}

A first instructive example is provided by catenoids: soap films, or minimal surfaces, formed between two coaxial parallel circular rings. In the elegant paper ``{\it In situ observation of a soap-film catenoid ---a simple educational physics experiment}'' by Ito and Sato~\cite{IS}, catenoids are experimentally produced in the laboratory and videotaped while the distance between the two circular wires is continuously increased. Besides catenoids, another critical point of the area functional always exists: the two flat disks spanned by the wires. For sufficiently small distances, two catenoids exist; among them, the one with a thicker neck is the absolute minimizer of the area (while the disks clearly have much larger area). As the wires separate, there is a distance $h_0$ at which the two states ---the thick-neck catenoid and the pair of disks--- have the same area. Just beyond this distance, the disks become the absolute minimizer, whereas the catenoid persists only as a local minimizer. Nevertheless, for an interval of distances $h>h_0$, the videotaped surface remains the catenoid, not the absolutely minimizing disks. This catenoid is thus a {\it stable minimal surface}, stability being understood as defined in the previous section.  

Throughout this process, the thinner-neck catenoid also exists, but it is always an unstable minimal surface (consistent with the general picture that a functional with two local minima should also possess a third, unstable, critical point). Eventually, at a larger critical distance $h_c$, the local minimizer (the thick-neck catenoid) and the unstable critical point (the thin one) merge, producing an inflection point. Beyond~$h_c$, the configuration of two disks is the only critical point. In the experiment~\cite{IS}, the unstable catenoid is captured in a photograph for a brief instant just before the distance $h_c$ is reached. Immediately afterwards, the thin neck collapses, and the soap film transforms into the two disks.  

The regularity theory of minimal surfaces has been the source of many major advances in the theory of elliptic and parabolic PDEs. In the 1960s, the Italian school proved that the Simons cone
\[
x_1^2+\ldots+x_m^2 \;=\; x_{m+1}^2+\ldots+x_{2m}^2 
\]
is an absolute minimizer of area (for its own boundary values in any ball of $\R^{2m}$) when $2m \geq 8$, while it is not even stable in dimensions $2$, $4$, and~$6$.\footnote{\label{expl_dim}This different behavior can be roughly explained by observing that the Jacobian for the area in spherical coordinates, $r^{2m-2}\,dr$, becomes smaller at the origin as the dimension $2m$ increases. For $2m=2$, the minimizer clearly avoids the origin: for the boundary values of the cone, it consists of two parallel lines, rather than the ``cross'' passing through the origin.}  
Thus, minimizing minimal surfaces of dimension $n$ may exhibit conical singularities when $n\geq 7$. At the same time, a sequence of outstanding contributions by various authors (culminated in J.~Simons' work) established that $n$-dimensional (absolute) minimizing minimal surfaces in $\R^{n+1}$ are always smooth when $n\leq 6$.  

Extending this regularity result to the larger class of {\it stable minimal surfaces} has been a long-standing open problem. Regularity was established for surfaces in $\R^3$ in the late 1970s, independently by Fischer-Colbrie and Schoen~\cite{Fischer-Schoen} and by do Carmo and Peng~\cite{doCarmo-Peng}. Much more recently, in the 2020s, the result has been extended to~$\R^4$ by Chodosh and Li~\cite{CL}, to $\R^5$ by Chodosh, Li, Minter, and Stryker~\cite{CLMS}, and to~$\R^6$ by Mazet~\cite{M}. The case of $\R^7$ (that is, $n=6$) remains open.

\section{Regularity of stable solutions to interior reaction problems}

The paper \cite{CFRS} takes on the analogue question, the regularity of {\it stable solutions}, for equations of the form $-\Delta u = f(u)$, where $\Delta$ is the Laplacian. They are called semilinear or reaction-diffusion elliptic equations and arise in many physical and biological situations. In the combustion problem described below, 
a similar phenomenon to that of  catenoids occurs. To describe it, let us start with the basic setting.

A much more detailed and self-contained exposition of the results of this section can be found in \cite{Cabre2022Holder}.

Consider the semilinear elliptic equation
\begin{equation}\label{eq:PDE}
-\Delta u=f(u) \qquad \text{in }\Omega\subset \R^n,
\end{equation}
where $\Omega$ is a bounded domain, $f\in C^1(\R)$, and  $u:\overline\Omega \subset \R^n\to \R$.
Notice that \eqref{eq:PDE} is the Euler-Lagrange equation of the functional
\begin{equation}\label{energy}
\mathcal E(u):=\int_{\Omega}\Bigl(\frac{1}{2}|\nabla u|^2-F(u)\Bigr)\,dx,
\end{equation}
where $F$ is a primitive of $f$. 
The second variation of $\mathcal E$ at $u$ is given by
$$
\frac{d^2}{d\varepsilon^2}\Big|_{\varepsilon=0}\mathcal E(u+\varepsilon\xi)
=\int_{\Omega}\Bigl(|\nabla \xi |^2-f'(u)\xi^2\Bigr)\,dx.
$$
We say that $u$ is a {\it stable} solution of equation \eqref{eq:PDE} in $\Omega$ if the second variation at $u$ is nonnegative, namely, if 
\begin{equation}\label{stabilityLip} 
\int_{\Omega} f'(u) \xi^2\,dx\leq \int_{\Omega} |\nabla \xi|^2\,dx  \quad \mbox{ for all }
\xi\in C^1(\overline\Omega) \text{ with } \xi_{|_{\partial\Omega}}\equiv 0.
\end{equation}
Note that the stability of $u$ is defined within the class of functions $u+\varepsilon\xi$ agreeing with~$u$ on the boundary $\partial\Omega$. Stability is equivalent to assuming the nonnegativeness of the first
Dirichlet eigenvalue in~$\Omega$ for the linearized operator of \eqref{eq:PDE} at~$u$, namely, $-\Delta-f'(u)$.

If a {\it local minimizer} of $\mathcal E$ exists (that is, a minimizer for small perturbations having same boundary values), then it will be a stable solution as well.

We are concerned with the regularity of stable solutions to \eqref{eq:PDE}, a question that has been investigated since the 1970s. It was motivated by the Gel'fand-type problems \eqref{eq:lambda} described below, for which $f(u)=e^u$ is a model case. Another typical nonlinearities to keep in mind are $f(u)=(1+u)^p$ with $p>1$. But we will see that our results will embrace very large classes of nonlinearities. Notice that for the above model nonlinearities, and more generally when $f$ satisfies \eqref{nonl-G} below, the energy functional \eqref{energy}, among functions with same boundary values  on~$\partial\Omega$ as any given function $u$, is unbounded below.\footnote{\label{Foot1}To see this,  consider $\mathcal E(u+tv)$ for any $v \in C^1_c(\Omega)$ with $v\geq0$ and $v\not\equiv 0$, and let $t\to +\infty$.}  Thus, the functional admits no absolute minimizer. However, we will see next that, in some relevant cases, it admits stable solutions.

Assume that $f:[0,+\infty)\to\R$ is $C^1$ and that
\begin{equation}\label{nonl-G}
f\text{ is nondecreasing, convex, and satisfies } f(0)>0 \text{ and } 
\lim_{t\to+\infty}\frac{f(t)}{t}=+\infty. 
\end{equation}
Given a constant $\lambda > 0$ consider the problem
\begin{equation}
\label{eq:lambda}
\left\{
\begin{array}{cl}
-\Delta u=\lambda f(u) & \text{in }\Omega\\
u>0 & \text{in }\Omega\\
u=0 & \text{on }\partial\Omega.\end{array}
\right.
\end{equation}
Note first that $-\Delta u=\lambda f(u)$ admits no trivial solution ($u\equiv 0$ is not a solution since $\lambda f(0)>0$). At the same time, the functional $\mathcal E$ is unbounded from below (see Footnote~\ref{Foot1}) and thus admits no absolute minimizer. However, \eqref{eq:lambda} admits a unique {\it stable} classical solution $u_\lambda$ for every $\lambda \in (0,\lambda^\star)$, where $\lambda^\star\in (0,+\infty)$ is a certain parameter; see e.g.\ the book \cite{Dup} by Dupaigne. On the other hand, no weak solution exists for $\lambda>\lambda^\star$ (under various definitions of weak solution). The functions $u_\lambda$, which are bounded in $\overline\Omega$ for $\lambda<\lambda^\star$ and hence regular, form an increasing sequence in $\lambda$. They converge, as $\lambda\uparrow\lambda^\star$, towards an $L^1(\Omega)$ distributional {\it stable} solution $u^\star$ of \eqref{eq:lambda} for $\lambda=\lambda^\star$. It is called {\it the extremal solution} of the problem and, depending on the nonlinearity~$f$ and the domain~$\Omega$, it may be either regular or singular.

When $f(u)=\lambda e^u$ (this is the so called {\it Gel'fand problem}), \eqref{nonl-G} is related to the thermal self-ignition of a chemically active mixture of gases in a container. The model was introduced by Frank-Kamenetskii in the 1930s but became popular within the mathematical community when Barenblatt wrote chapter~15 of the volume \cite{Gelf}, edited by Gel'fand in 1963.  Barenblatt chose to be $F(u)=\lambda e^u$, with $\lambda$ a positive constant, from Arrhenius law in chemical kinetics. Note that the temperature is kept at $u=0$ on the boundary $\partial\Omega$. Note the similarity of the structure of its solutions, as described above, with the situation for catenoids. Indeed, in the range $\lambda\in (0,\lambda^*)$ of parameters, the stable solution~$u_\lambda$ is not an absolute minimizer since $\mathcal E$ is unbounded by below. Instead, when $\lambda> \lambda^*$, no solution exists ---in the same way that catenoids did not exist for distances $h$ between the wires larger than $h_c$.

To better understand the problem, let us consider the nonlinear heat equation 
\begin{equation}
\label{eq:heat}
v_t-\Delta v=\lambda f(v), 
\end{equation}
where $v=v(x,t)$ and $t$ is time. Now, a stable solution $u=u(x)$ of \eqref{eq:PDE} can be understood as a stationary solution of \eqref{eq:heat} which is stable in the sense of Lyapunov ---note that a simple computation shows that the action functional $\mathcal E(v(\cdot,t))$ is non-increasing in the time $t$. The problem is nonlinear due to the sources of heat, $f(v(x,t))$ or $f(u(x))$: the production of heat depends nonlinearly on the actual temperature. As described in \cite{GM,GMN}, equation~\eqref{eq:heat}
describes the evolution of an initially uniform temperature $v(\cdot,0)\equiv 0$ which diffuses in space and increases in the container due to the heat release given by the reaction term $\lambda f(v)$ ---note that in Gel'fand's problem the initial heat source $\lambda f(0)= \lambda e^0=\lambda$ is already positive. The parameters $\lambda$ for which there exists a stable solution  of \eqref{eq:PDE} correspond to ignition failure (the reactive component undergoes partial oxidation and results in establishing a stationary temperature profile equal to the stable solution). Instead, $\lambda>\lambda^*$ (when there exists no stationary solution) means successful auto-ignition in the combustion process.

\subsection{The interior regularity result}\label{subsec-interior-results}
\
\smallskip

In our setting, as for the Simons cone in minimal surfaces treated in Section~2, there is an explicit singular stable solution (with finite energy) in large dimensions. Indeed, for $n\geq 3$,
\begin{equation}\label{logsoln}
u=\log\frac{1}{|x|^2} \quad \text{solves \eqref{eq:PDE} weakly, with } f(u)=2(n-2)e^u \text{ and } \Omega=B_1,
\end{equation}
and vanishes on $\partial B_1$. Using Hardy's inequality,\footnote{See Proposition 1.20 in \cite{CP}, for instance, for its very simple proof.} we see that the linearized operator $-\Delta-2(n-2)e^u=-\Delta-2(n-2)|x|^{-2}$ is nonnegative when $2(n-2)\leq (n-2)^2/4$, i.e., when $n\ge 10$. Thus, there exist singular $W^{1,2}(B_1)$ stable solutions to equations of the form \eqref{eq:PDE} whenever $n\geq 10$. From our results we will see that this does not occur for $n\leq 9$.

In the 1970s the seminal paper of Crandall and Rabinowitz~\cite{CR} established the boundedness of~$u^\star$ when $n \leq 9$ and $f$ is either the exponential $f(u)=e^u$ or a power $f(u)=(1+u)^p$ with $p>1$. In the mid-nineties Brezis asked for a similar result in the larger class of nonlinearities \eqref{nonl-G} (see, e.g., \cite{Brezis}), a question that gave rise to the following works, among others, cited in chronological order.

Nedev~\cite{Ned00} proved the regularity of $u^\star$ for $n\leq 3$. In the radial case,
Cabr\'e and Capella~\cite{CC} established the boundedness
of every stable radial $W^{1,2}(B_1)$  solution to \eqref{eq:PDE} in the unit ball
whenever $n\leq 9$, for every nonlinearity $f$.
Back to the general nonradial case, an interior $L^\infty$ a priori bound for stable regular solutions holding  for all nonlinearities was established for $n\leq 4$ by Cabr\'e~\cite{C10} in 2010. The author gave a different proof of this result more recently, in~\cite{C19}. This led to the regularity of $u^\star$ up to the boundary of convex domains when $n\leq 4$. The convexity assumption of the domain was nicely removed by Villegas in~\cite{Vil13}.

The optimal dimension $n\leq 9$ was finally reached in 2020, by Cabr\'e, Figalli, Ros-Oton, and Serra~\cite{CFRS}, assuming $f\ge 0$ for interior regularity, and $f\ge 0$, $f'\ge 0$, and $f''\ge 0$ for boundary regularity.\footnote{An heuristic explanation of why a maximal dimension for regularity exists, in this case being~9, is the one given previously in Footnote~\ref{expl_dim}.}
The following is the interior regularity result of \cite{CFRS}. It provides an interior H\"older a priori bound for stable solutions  when $n\leq 9$. It also establishes, in every dimension, an interior $W^{1,2}$ estimate (indeed a better $W^{1,2+\gamma}$ bound). All quantities are controlled in terms of only the $L^1$ norm of the solution. Regarding the nonlinearity, all what we need is $f\geq 0$. Here, and throughout the paper, by ``dimensional constant'' we mean a constant that depends only on $n$. 

\begin{theorem}[\cite{CFRS}, Theorem 1.2]\label{thm:0}
Let $u\in C^\infty(\overline B_1)$ be a stable solution of $-\Delta u=f(u)$ in $B_1\subset \R^n$, for some nonnegative function $f\in C^{1}(\R)$.

Then,
\begin{equation}
\label{eq:W12 L1 int}
\|\nabla u\|_{L^{2+\gamma}({B}_{1/2})} \le C\|u\|_{L^1(B_1)} \quad \text{ for every } n\geq 1
\end{equation}
for some dimensional constants $\gamma>0$ and $C$.
In addition,
\begin{equation}
\label{eq:Ca L1 int}
 \|u\| _{C^\alpha(\overline{B}_{1/2})}\leq C\|u\|_{L^1(B_1)} \quad \text{ if } n \leq 9,
 \end{equation}
where $\alpha>0$ and $C$ are dimensional constants.
\end{theorem}

The theorem is stated as an a priori bound for smooth stable solutions. From it one may deduce, by an approximation method, interior regularity for $W^{1,2}$ stable weak solutions to \eqref{eq:PDE} whenever $f$ is nonnegative, nondecreasing, and convex. The same fact holds for the extremal solution~$u^\star$ of~\eqref{eq:lambda} when $f$ satisfies~\eqref{nonl-G}.

For $n\geq 10$ (and still $f\geq 0$), \cite{CFRS} proved that stable solutions belong, in the interior of the domain, to certain Morrey spaces. This was an almost optimal result that has been improved by Peng, Zhang, and Zhou \cite{PZZ1}, reaching optimality. For instance, for $n=10$, \cite{PZZ1} shows that stable solutions are BMO in the interior, a result also found, independently, by Figalli and Mayboroda (personal communication).

\subsection{The boundary result}\label{subsec-bdry-results}
\
\smallskip

Let us now turn into boundary regularity.
It is well known that the interior bound of Theorem~\ref{thm:0} and the moving planes method yield an $L^\infty$ bound up to the boundary when the domain~$\Omega$ is convex and we deal with vanishing Dirichlet boundary conditions; see \cite{C10} and Corollary~1.4 of \cite{CFRS}. Since this procedure does not work for nonconvex domains, \cite{CFRS} undertook the study of boundary regularity in more general domains and proved the following result.

\begin{theorem}[\cite{CFRS}, Theorem 1.5]
\label{thm:7}
Let  $\Omega\subset \R^n$ be a bounded domain of class $C^3$ and let $u\in C^0(\overline\Omega) \cap C^2(\Omega)$ be a stable solution of 
\begin{equation}\label{BVP}
\left\{
\begin{array}{cl}
-\Delta u=f(u) & \text{in }\Omega\\
u=0 & \text{on }\partial\Omega.
\end{array}
\right.
\end{equation}
Assume that $f: \R \to \R$ is nonnegative, nondecreasing, and convex.

Then,
\begin{equation}\label{eq:w12}
\|\nabla u\|_{L^{2+\gamma}(\Omega)}\le C_\Omega \,\|u\|_{L^1(\Omega)} \quad \text{ for every } n\geq 1,
\end{equation}
and
\begin{equation*}
 \|u\| _{C^\alpha(\overline\Omega)} \leq C_\Omega\, \|u\|_{L^1(\Omega)} \quad \text{ for } n\leq 9,
\end{equation*}
where $\gamma>0$ and $\alpha>0$ are dimensional constants, while $C_\Omega$ depends only on $\Omega$.
\end{theorem}

Erneta~\cite{Ern1, Ern2, Ern3} has improved the previous result to allow for~$C^{1,1}$ domains.
 
From the a priori estimates of the theorem, one deduces the corresponding regularity results for the extremal solution $u^\star$ to \eqref{eq:lambda}. 

The $W^{1,2}$ regularity result \eqref{eq:w12} solved~\cite[Open problem 1]{BV}, a question posed by Brezis and V\'azquez. Previously, it had been proven, under the stronger hypotheses~\eqref{nonl-G} on the nonlinearity, for $n\leq 5$ by Nedev~\cite{Ned00} and for $n=6$ by Villegas~\cite{Vil13}.

\subsection{The three ingredients for the interior estimate}\label{subsec-interior}
\
\smallskip

We first note that, by approximation, the stability inequality \eqref{stabilityLip},
$$
\int_{\Omega} f'(u) \xi^2\,dx\leq \int_{\Omega} |\nabla \xi|^2\,dx,
$$ 
holds not only for $C^1(\overline\Omega)$ functions $\xi$ vanishing on $\partial\Omega$, but also for all Lipschitz functions $\xi$ in $\overline\Omega$ which vanish on $\partial\Omega$. Hence, we can consider a test function of the form $\xi=\mathbf{c}\eta$,
where $\mathbf{c}\in W^{2,\infty}(\Omega)$, $\eta$ is a Lipschitz function in $\overline\Omega$, and $\mathbf{c}\eta$ vanishes on~$\partial\Omega$. 
Then, 
since
\begin{equation}\label{eq:parts}
\int_{\Omega} |\nabla \xi|^2\,dx = \int_{\Omega} \bigl(  |\nabla \mathbf{c}|^2 \eta^2 + \mathbf{c}\nabla\mathbf{c} \nabla \eta^2+  \mathbf{c}^2 |\nabla \eta|^2\bigr) \, dx,
\end{equation}
integrating by parts the second term of the last integral we see that 
\begin{equation}\label{eq:07}
\int_{\Omega} \bigl( \Delta \mathbf{c}+f'(u)\mathbf{c}\bigr) \mathbf{c}\,\eta^{2}\, dx \leq 
\int_{\Omega} \mathbf{c}^{2}|\nabla \eta|^{2} \, dx.
\end{equation}

The key point here is that, if we take $\eta$ vanishing on the boundary, then we have the freedom to make choices for the function $\mathbf{c}$ having arbitrary values on $\partial\Omega$.
Note also that expression \eqref{eq:07} brings the linearized operator $\Delta+f'(u)$, acting on $\mathbf{c}$, into play. Later there will be two crucial different choices for the function~$\mathbf{c}$: \eqref{test1} and \eqref{test2} below.

The proof of H\"older regularity is based on putting together three estimates. In the interior case, the first one, and where the restriction $n\leq 9$ appears, is the following. Here $f\geq 0$ is not needed. 
Throughout the paper we will use the notation
\begin{equation}\label{notation r}
r=|x| \qquad\text{and}\qquad u_r=\frac{x}{|x|}\cdot \nabla u.
\end{equation}

\begin{proposition}[\cite{CFRS}, Lemma 2.1] \label{prop:1}
Let $u\in C^\infty (\overline B_1)$ be a stable solution of $-\Delta u=f(u)$ in $B_1\subset \R^n$, for some function $f\in C^{1}(\R)$.

If $3 \leq n\leq 9$, then 
\begin{equation}\label{ineq1}
\int_{B_\rho}r^{2-n} u_{r}^2\,dx \leq C\rho^{2-n}\int_{B_{3\rho/2}\setminus B_\rho}|\nabla u|^2\,dx
\end{equation}
for all $\rho< 2/3$, where $C$ is a dimensional constant.
\end{proposition}

This result is proven using \eqref{eq:07} with 
\begin{equation}\label{test1}
\mathbf{c}= x\cdot\nabla u = ru_r \quad\text{ and } \quad \eta=r^{(2-n)/2}\zeta,
\end{equation}
where $\zeta$ is a cut-off function. Recall that $r=|x|$.

Note that the proposition requires $n\geq 3$. However, adding superfluous independent variables to the solution, as in \cite{CFRS}, one sees that it can be used to establish Theorem~\ref{thm:0} also in dimensions one and two.

If the right-hand side of estimate \eqref{ineq1} had $u_r^2$ as integrand, instead of $|\nabla u|^2$, then H\"older regularity of $u$ would be a simple consequence of it, by the ``hole filling technique''; see \cite{CFRS}. To overcome this difficulty we need two more ingredients.

The first one uses the choice 
\begin{equation}\label{test2}
\mathbf{c} = |\nabla u|
\end{equation}
in \eqref{eq:07} and leads to certain $L^2$ and $L^1$ second derivative estimates. From them, we will deduce the $L^{2}$ estimate for $\nabla u$ of the following proposition, which holds in all dimensions.

\begin{proposition}[\cite{CFRS}, Proposition 2.5]\label{prop:2}  
Let $u\in C^\infty (\overline B_1)$ be a stable solution of $-\Delta u=f(u)$ in $B_1\subset \R^n$, for some nonnegative function $f\in C^{1}(\R)$.

Then,
\begin{equation}\label{ineq2}
\|\nabla u\|_{L^{2}(B_{1/2})}  \le C \|u\|_{L^{1}(B_{1})}
\end{equation}
for some dimensional constant $C$.
\end{proposition}

Proposition \ref{prop:2} provides the second step towards the H\"older regularity proof. Indeed, one starts from estimate \eqref{ineq1} of Proposition~\ref{prop:1} with, for instance,  $\rho=5/8$. Now, one covers the annulus in the right-hand side of \eqref{ineq1} by a dimensional number of balls with appropriate small radius. Next, in each of these balls one uses \eqref{ineq2}, properly rescaled, with $u$ replaced by the stable solution $u-t$, a solution of $-\Delta v= f(v+t)$, where $t$ is a constant to be chosen later. Adding the right-hand sides of  \eqref{ineq2} in all the balls, one concludes that the integral of $r^{2-n} u_{r}^2$ in $B_{1/2}$ can be controlled by the $L^1$ norm of $u-t$ in the annulus $B_1\setminus B_{1/2}$.

To control this last quantity, we will use the final third ingredient, namely, estimate \eqref{new-int-ur} in an annulus of Proposition~\ref{prop:3} below, first established by Cabr\' e in~\cite{C22quant}. The paper \cite{CFRS} used instead another related result which did not have such a ``clean'' statement. In addition, it was proved there by a contradiction-compactness proof. Instead, in~\cite{C22quant}, Proposition~\ref{prop:3} is proved in a quantitative manner.

\begin{proposition}[\cite{C22quant}, Theorem 1.4] \label{prop:3}
Let $u\in C^\infty (\overline B_1)$ be superharmonic in the ball $B_1\subset \R^n$. 

Then, 
\begin{equation}\label{new-int-ur} 
\big\Vert u-\textstyle{\ave_{B_{1}\setminus B_{1/2}}}u\big\Vert_{L^{1}({B_{1}\setminus B_{1/2})}}\leq C\Vert u_{r}\Vert_{L^{1}({B_{1}\setminus B_{1/2})}}
\end{equation}
and
\begin{equation}\label{new-int-ur-ball} 
\big\Vert u-\textstyle{\ave_{B_{1}}}u\big\Vert_{L^{1}({B_{1})}}\leq C\Vert u_{r}\Vert_{L^{1}({B_{1})}}
\end{equation}
for some dimensional constant $C$.
\end{proposition}

Note now that the three propositions, used one after the other as described above, lead to \eqref{ineq1} with $\nabla u$ replaced by $u_r$ on its right-hand side. Then the ``hole filling technique" easily leads to H\"older regularity.

\subsection{On the radial Poincar\'e inequality}

One could wonder about the validity of \eqref{new-int-ur} and  \eqref{new-int-ur-ball} for all functions~$u$, and not only superharmonic ones. For $n\geq2$, both inequalities would be false. Indeed, this can be seen with the family of smooth functions 
$$
u^\ep(x):=\varphi(|x|/\ep)\, x_n/|x|, 
$$
where $\ep\in (0,1)$, $\varphi$ is smooth and nonnegative, vanishes in $(0,1/4)$, and is identically~1 in $(1/2,\infty)$. Then, $u^\ep_r$ has support in $B_{\ep/2}\setminus B_{\ep/4}$ and its $L^1(B_1)$ norm is bounded above by $C\ep^{n-1}$. Instead, to control the $L^1(B_1\setminus B_{1/2})$ norm of $u^\ep-\textstyle{\ave_{B_1\setminus B_{1/2}}} u$ (or, similarly, of $u^\ep-\textstyle{\ave_{B_1}} u$) by below, if $\textstyle{\ave_{B_1\setminus B_{1/2}}} u\leq0$ (the case $\textstyle{\ave_{B_1\setminus B_{1/2}}} u\geq 0$ is treated similarly) we may integrate only in $\{1/2<|x|<1, x_n>0\}$. In this set $|u^\ep-\textstyle{\ave_{B_1\setminus B_{1/2}}} u|=u^\ep -\textstyle{\ave_{B_1\setminus B_{1/2}}} u\geq u^\ep=x_n/|x|$, and thus the $L^1(B_1\setminus B_{1/2})$ norm is larger than a positive constant independent of $\ep$.

Instead, it is easy to check (this will be done in \cite{CSa}) that the radial Poincar\'e inequalities \eqref{new-int-ur} and \eqref{new-int-ur-ball} hold among all functions that vanish on $\partial B_1$ (in this case, there is no need to subtract a constant to the function $u$ for the inequalities to hold). 

In the forthcoming paper \cite{CSa}, Cabr\' e and Saari determine whether the radial Poincar\'e inequalities above hold in $L^p$ instead of $L^1$, depending on the dimension and on~$p$.

\begin{theorem}[\cite{CSa}]\label{Lp-version}
Given $p\in [1,\infty)$, inequality \eqref{new-int-ur} in the annulus holds, for all superharmonic functions, with $L^1$ replaced by $L^p$  in both its left and right-hand sides if and only~if
\begin{itemize}
\item \text{either } $n\leq 3 \text{ and } 1\leq p<\infty$,
\item \text{or } $n\geq 4 \text{ and } 1\leq p <(n-1)/(n-3)$.
\end{itemize}

In addition, in the range of exponents where the inequality does not hold $($that is, $p\geq (n-1)/(n-3))$, there is a sequence of smooth superharmonic functions which provide a counterexample to the validity of the inequality, and which converge to the function 
$$
u(x)= cr^a + G(\sigma),
$$ 
where $x=r\sigma$ with $\sigma \in S^{n-1}$, $c$ and $a<0$ are appropriate constants, and $G$ is the Green's function of the spherical Laplacian.
\end{theorem}

We have two different proofs of the result. The first one is quite elementary: it decomposes the Laplacian in its radial and spherical parts, and later uses integration by parts and a duality argument together with elliptic regularity on the sphere. The second proof uses the Fourier decomposition in spherical harmonics.

\subsection{Related equations. Open problems}\label{subsec-open-related}
\

\smallskip
3.5.1. The equation $(-\Delta)^su=f(u)$, with $0<s<1$,  involving the fractional Laplacian, will appear in next sections in relation with boundary reaction problems. For it, the optimal dimension for regularity of stable solutions has only been reached when $f(u)=\lambda e^u$. This was accomplished in a remarkable paper by Ros-Oton~\cite{R-O}. The dimension found in that paper, which depends on~$s\in (0,1)$, is expected to also give the correct range of dimensions for other nonlinearities. This is however an open problem, even in the radial case. The most recent progress is the paper~\cite{CS-P} by Cabr\' e and Sanz-Perela in the general nonradial case, where the techniques of \cite{CFRS} are extended for the half-Laplacian ($s=1/2$) to prove interior regularity in dimensions $n\leq 4$ for all nonnegative convex nonlinearities ---here, with $s=1/2$, one expects $n\leq 8$ to be the optimal result.

\smallskip
3.5.2. The interior result of \cite{CFRS} has been extended to the $p$-Laplacian by Cabr\'e, Miraglio, and Sanch\'on~\cite{CMS}. This work reaches the optimal dimension for interior regularity of stable solutions when $p>2$. It remains an open problem to accomplish the expected optimal dimension for $1<p<2$.

\smallskip
3.5.3. Bruera and Cabr\'e~\cite{BC} have extended some of the techniques of \cite{CFRS} to stable solutions of MEMS (micro-electro-mechanical systems) problems to reach a regularity result in the optimal range $n\leq 6$. Here one has the equation $-\Delta u =\lambda f(u)$ for positive nonlinearities $f:[0,1)\to\R$ which blow-up at $u=1$ in such a way that $\int_0^1 f=+\infty$. In this new setting, a solution is said to be regular when $\Vert u\Vert_{L^\infty(\Omega)}< 1$. The result requires, however, a Crandall-Rabinowitz type assumption on the nonlinearity, more precisely, that $\liminf_{t\uparrow 1} \left(ff''/(f')^2\right) (t) >1$ ---or, equivalently, that $f^{-a}$ is concave near~1 for some $a>0$. It is an open question to decide if such Crandall-Rabinowitz condition can be removed.

\smallskip
3.5.4.  After \cite{CFRS}, an interesting result by Peng, Zhang, and Zhou~\cite{PZZ2} proves that, for $n=5$ and the Laplacian, the hypothesis $f\geq 0$ is not needed to conclude interior H\"older regularity. The H\"older control, however, is given by the $W^{1,2}$ norm of the solution instead of the $L^1$ norm. This establishes in dimension~5 the previous analogue result of Cabr\'e~\cite{C10,C19} for $n\leq 4$. In the general nonradial case, it is not known if an interior $L^\infty$ estimate could hold for $6\leq n\leq 9$ without the hypothesis $f\geq 0$.

\smallskip
3.5.5. Our boundary regularity results require $f\geq 0$, $f'\geq 0$, and $f''\geq 0$. Both the proof in \cite{CFRS} and the one from \cite{C22quant}, which are completely different, turn out to need these three conditions. It is not known if these hypotheses could be relaxed.

\section{Stable solutions to boundary reaction problems}

We now turn into reaction problems where the boundary value is not prescribed. That is, minimizers and stable solutions are considered in the larger space $H^1(\Omega)$ instead of (in the case, for instance, of zero Dirichlet boundary values) the smaller space $H^1_0(\Omega)$. Consequently, critical points satisfy a Neumann boundary condition (a value for the normal derivative $\partial_{\nu} u $) instead of the Dirichlet one.

\subsection{Interior reactions under zero Neumann boundary conditions}

The following is a  celebrated result from the late 1970s, proved by Casten and Holland, and independently by Matano. It establishes the nonexistence of nonconstant stable solutions
of the Neumann problem in convex domains of $\R^n$.

\begin{theorem}[\cite{CastenHolland1978,Matano1979}]\label{thm:Casten-Holland}
Let $\Omega\subset\R^n$ be a bounded, smooth, convex domain, and let $f\in C^1(\R)$. 
Assume $u\in C^2(\overline\Omega)$ is a solution of the Neumann problem
\begin{equation}\label{BVP-Neumann}
\left\{
\begin{array}{cl}
-\Delta u = f(u) & \text{in }\Omega \\[6pt]
\partial_{\nu} u = 0 & \text{on }\partial\Omega
\end{array}
\right.
\end{equation}
and that $u$ is stable, in the sense that
\[
\int_{\Omega} f'(u)\,\xi^2\,dx \le \int_{\Omega} |\nabla\xi|^2\,dx
\qquad\text{for all }\xi\in C^1(\overline\Omega).
\]

Then, $u$ is constant.
\end{theorem}

The proof of the result is simple. To understand its statement, it is interesting to look at {\it the Allen-Cahn equation}
\begin{equation}\label{eq:pde_int}
\left\{
\begin{array}{cl}
- \Delta u = \dfrac{1}{\varepsilon^2} (u-u^3) & \text{in }\Omega \\[2mm]
\partial_{\nu} u = 0 & \text{on }\partial\Omega,
\end{array}
\right.
\end{equation}
which is the first variation of the functional 
\begin{equation}\label{int_ener}
\hat E_\varepsilon(u)=\int_\Omega\left( \frac12|\nabla u|^2 + \frac1{4\varepsilon^2} (1-u^2)^2\right)\, dx.
\end{equation}
Here the potential energy corresponds to a double-well potential. Critical points will satisfy $-1\leq u\leq 1$.

A well-known $\Gamma$-convergence result, \cite{Modica1987,KohnSternberg1989},  states that the rescaled energy $\varepsilon\hat E_\varepsilon$ converges, as $\varepsilon\to 0$,  to a multiple of the perimeter of the transition set where characteristic functions $u$ change from $-1$ to $1$. More precisely,
the limit energy is 
\begin{equation}\label{eq:energy_hat}
\hat{E}(u) =
\begin{cases}
c\,\mathrm{Per}_\Omega(\{ u = 1\}) & \text{if } u(x) = \pm 1 \text{ a.e. in }\Omega \\[1mm]
+\infty & \text{otherwise}.
\end{cases}
\end{equation}
To give some light into Theorem~\ref{thm:Casten-Holland}, note that given a strictly convex domain (to simplify in the plane $\R^2$) and any straight line segment (being the limiting transition layer in $\Omega$ where the solution goes from $-1$ to 1 ---we take it to be a straight line since it must minimize the perimeter or length functional) connecting two points on~$\partial\Omega$, it is always possible to find a shorter straight line connecting two other near-by points on the boundary, by the convexity of $\Omega$. This gives nonexistence of local minimizers for the limiting problem, suggesting that the same could be true for the Allen-Cahn problem for $\epsilon$ small. This fact is indeed true by Theorem~\ref{thm:Casten-Holland}.

Instead in non-convex domains, for example in dumbbell-shaped ones, it is possible to construct nonconstant stable solutions of \eqref{eq:pde_int} near a local minimizer of 
the limiting functional, as done by Kohn and Sternberg~\cite{KohnSternberg1989}. Note that in many dumbbell-shaped domains there will be a unique straight line segment connecting two opposite points in the thin neck with such straight segment being a strict local minimizer of the length functional.

\subsection{Boundary reaction problems}

Let us now turn into the case of {\it boundary reactions}. Here, the potential energy is computed on the boundary $\partial\Omega$ instead of~$\Omega$. Let us stay in dimension 2, so that $\Omega$~is a smooth, bounded, simply connected domain in $\R^2$. Our functional is now
\begin{equation}\label{eq:ener}
E_{\varepsilon}(u) =  \int_{\Omega} \frac{1}{2}|\nabla u|^2 \, dx +
\int_{\partial \Omega} \frac{1}{4\varepsilon}(1-u^2)^2 \, d{\ell},
\end{equation} 
instead of \eqref{int_ener}. Its first variation is the boundary reaction problem
\begin{equation}\label{eq:pde}
\left\{
\begin{array}{cl}
\Delta u = 0 & \text{in }\Omega \\[6pt]
\partial_{\nu} u = \dfrac{1}{\varepsilon}(u-u^3) & \text{on }\partial\Omega,
\end{array}
\right.
\end{equation}
where $\varepsilon>0$.

For this problem, the limit $\varepsilon\to 0$ is quite different than in the previous interior case. Along the boundary, the transition sets (to go from $-1$ to 1) converge to \textit{points},
and the relevant $\Gamma$-limit only counts the number of points, which is discrete and hence does not have nontrivial local minima. Hence, the limiting problem gives no information about the existence, or not, of local minimizers. This is in contrast with the interior reaction case where existence, respectively nonexistence, was obvious, for the limiting problem, in dumbbell-shaped domains, respectively in strictly convex domains.

In this respect, C\'onsul~\cite{Consul1996} proved the nonexistence of nonconstant stable solutions to \eqref{eq:pde} when the domain is a ball of $\R^n$, even when the reaction term $(u-u^3)/\varepsilon$ is replaced by any nonlinearity $f(u)$. On the other hand, the analogue of the interior dumbbell result was established in 1999 by C\'onsul and Sol\`a-Morales~\cite{ConsulSola-Mor1999}. They proved the existence of nonconstant stable solutions to \eqref{eq:pde} when $\Omega$ are appropriate dumbbell domains.
  
Since then, the validity of the analogue result to Theorem~\ref{thm:Casten-Holland} for boundary reactions remained open. This question was raised by Sol\`a-Morales in the mid 1990s.  It has been recently solved by Cabr\'e, C\'onsul, and Kurzke~\cite{CCK}, who prove that such analogue is not true, even in some simple domains in $\R^2$. Its main result is the following.

\begin{theorem}[\cite{CCK}]\label{thm:T0}
Let  $\Omega$ be the square $(-1,1)\times (-1, 1)$ in $\R^2$, or some smooth strictly convex domains approximating the square. 

Then, for $\varepsilon$ sufficiently small, there exists a nonconstant stable solution to problem \eqref{eq:pde}. In addition, such solutions converge on $\partial\Omega$, as $\varepsilon\to 0$ and in the case of the square, to the function which is $-1$ on the left of the centered points $(0,-1)$ and $(0,1)$, and is $1$ on  the part of $\partial \Omega$ in their right.
\end{theorem}

This result originated from numerical computations of such solutions done by C\'onsul and Jorba~\cite{Consul2005}. They were found through a bifurcation-continuation method. Their very fine finite element method allowed for checking, numerically, that the solutions making the transition near the two centered points of the square were stable for $\varepsilon$ small enough.

Our proof of Theorem~\ref{thm:T0} is quite elaborate and requires the development of a Ginzburg-Landau theory as the one in Bethuel-Brezis-Helein~\cite{BethuelBrezisHelein1994}. The functions involved are now scalar valued (in contrast with the complex valued in~\cite{BethuelBrezisHelein1994}) and take $S^0=\{ -1,1\}$ values on the boundary $\partial\Omega$. 

Let us explain this in more detail. Given two points $p$ and $q$ on $\partial\Omega$, $p\neq q$, let $I_{p,q}$ be the connected component of $\partial\Omega\setminus\{p,q\}$ corresponding to going from $p$ to $q$ counterclockwise. We define $\chi^{p,q}$ as the characteristic  function on the boundary given by
\begin{equation}\label{eq:chif}
\chi^{p,q} =
\begin{cases}
-1 & \text{ on } I_{p,q}  \\
1 & \text{ on } \partial\Omega\setminus  I_{p,q}.
\end{cases}
\end{equation}
Note that this function
$\chi^{p,q}$ minimizes the boundary term in the energy \eqref{eq:ener} since $(1-(\chi^{p,q})^2)^2=0$ on all $\partial\Omega$. Hence, it could be natural
to guess that a minimizer of the total energy would be obtained by the harmonic extension to
$\Omega$ of $\chi^{p,q}$, that is, $u_0^{p,q}$ satisfying
\begin{equation}\label{hext}
\left\{ \begin{array}{rcll} \Delta u_0^{p,q} & = & 0 &\mbox{ in } \Omega \, , \\
u_0^{p,q}& = & \chi ^{p,q}  & \mbox{ on } \partial\Omega\, .
\end{array} \right.
\end{equation}
Notice that since $\chi^{p,q}\not\in H^{1/2}(\partial\Omega)$, we have $u_0^{p,q}\not\in H^1(\Omega)$. So, it
can not be a minimizer of the energy. Indeed, the Dirichlet energy of $u_0^{p,q}$ is infinite, but we can still study its dependence on the choice of points $p,q$ by renormalizing appropriately.
More precisely, as in the complex valued Ginzburg-Landau theory, we can remove small balls around $p$ and $q$ and study the behavior of the Dirichlet energy on $\Omega_\rho=\Omega\setminus (B_\rho(p) \cup
B_\rho(q))$ as $\rho\to 0$. In~\cite{CCK} we prove that
\begin{equation}\label{eq_ener}
\frac12\int_{\Omega_{\rho}} |\nabla u_0^{p,q}|^2 \, dx= \frac 4\pi  \log \frac{1}{\rho} +
W_{\Omega}(p,q) + O(\rho), 
\end{equation}
where $W_\Omega$ is called the {\it renormalized energy} of the problem and can be computed by means of the conformal map from $\Omega$ to the half-space or from Green's functions of~$\Omega$. Indeed, we have that
\begin{align}\label{renorm}
W_{\Omega}(p,q)&=-\frac{2}{\pi} \log \frac{\partial ^2 G^D}{\partial \nu_p
\partial \nu_q}(p,q) \\
&=2\left( R^N(p,q)+R^N(q,q)-2 G^N(p,q)\right)\\
&= {\frac{4}{\pi} \log |p-q| + 2 \left( R^N(p,p)+R^N(q,q)-2R^N(p,q)\right) },
\end{align}
where $G^D(p,q)$ is the Dirichlet Green's function of
$\Omega$, $G^N(p,q)$ is the Neumann Green's function of $\Omega$, and $R^N(p,q)$ is the regular
part of the last one.

The proof of Theorem \ref{thm:T0} also requires to prove lower and upper bounds similar to \eqref{eq_ener}, but with $\rho=\varepsilon$ and $u_0^{p,q}$ replaced by minimizers $w_\e$ of the functional $E_\e$ among functions having boundary values in a small {\it closed} $L^2(\partial\Omega)$ neighborhood of~$\chi^{p,q}$. More precisely, the lower bound is
\begin{equation}\label{lowerb}
E_\e(w_\e) \ge \frac4\pi \log\frac1\e + W_\Omega(p,q) +2C_f -o_\e(1),
\end{equation}
while the upper bound is the same with $-o_\e(1)$ replaced by $+o_\e(1)$.

Assume now that the domain $\Omega$ admits a couple of points $p$ and $q$ on its boundary such that $W_\Omega:\partial\Omega\times\partial\Omega\longrightarrow \R$ has at $(p,q)$ a strict isolated local minimizer. Then with the aid of the above lower and upper bounds, one can prove that the energy \eqref{eq:ener} admits, for $\e $ small, a local minimizer that will converge, as $\e\to 0$, to the characteristic function $\chi^{p,q}$. Such local minimizers are the claimed nonconstant stable solutions to \eqref{eq:pde} in the theorem.

Thus, the remaining (simple) task in \cite{CCK} consists in showing that the renormalized energy \eqref{renorm}, with $\Omega$ being the square, has a strict local isolated minimizer at the couple of opposite centered points $\left( (0,-1), (0,1) \right)$.

To establish \eqref{lowerb}, and also the upper bound, one blows-up problem \eqref{eq:pde} near the point $p$, obtaining the boundary reaction problem 
\begin{equation}\label{eq:halfhalf}
\left\{
\begin{array}{cl}
\Delta u = 0 & \text{in }\R^2_+:=\{x\in\R, y>0\} \\[6pt]
\partial_{\nu} u = f(u) & \text{on }\partial\R^2_+
\end{array}
\right.
\end{equation}
in the half-plane, with $f(u)=u-u^3$. One then uses the results on stable and on increasing solutions to problem \eqref{eq:halfhalf} obtained in 2005 by Cabr\'e and Sol\`a-Morales~\cite{CabreSolaMorales2005}. In that paper it was very useful to know that the function
\[
\phi(x,y)=\dfrac{2}{\pi}\arctan \frac{x}{y+1}
\]
is a solution to \eqref{eq:halfhalf} when $f(u)=\sin(\pi u)$ ---this is called the Peierls-Nabarro problem on crystal dislocations.

Note also that problem \eqref{eq:halfhalf} is equivalent, by the so called Caffarelli-Silvestre extension technique and with $v=u(\cdot,0)$, to the equation
\begin{equation}\label{half-lapl}
(-\Delta)^{1/2} v = f(v) \quad\text{in } \R
\end{equation}
for the half-Laplacian, an operator that we treat in the next and last section.

\section{Reaction equations for the fractional Laplacian on the whole real line. Periodic solutions}

It is well known that bounded solutions to the semilinear second order ODE $-u''=f(u)$ in all of $\R$ are ---up to a multiplicative factor~$\pm 1$ and up to translations--- either increasing in $\R$, or even with respect to $0$ and decreasing in $(0,+\infty)$, or periodic. This follows immediately from the fact that $u$ must be even with respect to any of its critical points\footnote{This is a consequence of the uniqueness theorem for ODEs: if $u'$ vanishes at $x=0$, say, then $v(x):=u(-x)$ solves the same equation and has same initial position and derivative at $x=0$ as $u$.} after considering the cases when $u'$ vanishes at none, only one, or at least two points. 

It is then natural to ask whether this classification  also holds true for the fractional Laplacian, that is, if, for every 
$0<s<1$, any bounded solution $u$ to 
\begin{equation}
 \label{eq:frac Lapl intro}
  (-\Delta)^s u=f(u) \quad\text{ in }\R
\end{equation}
belongs to one of the above mentioned three categories. Recall that
\begin{equation*}\label{defi laplacian_}
(-\Delta)^s u(x):=c_s\, \text{P.V.}\!\int_\R 
\frac{u(x)-u(y)}{|x-y|^{1+2s}}\, dy,\qquad
c_s:=\frac{s4^s\Gamma(1/2+s)}{\sqrt\pi\,\Gamma(1-s)}.
\end{equation*} 
Due to the lack of an analogue uniqueness result for  \eqref{eq:frac Lapl intro}, the aforementioned argument cannot be carried out in the nonlocal case.  

This question, which is still open, is due to Sol\`a-Morales, who posed it  for the case $s=1/2$, or equivalently for problem \eqref{eq:halfhalf}. The conjecture is only known to be true for the equation
\begin{equation}\label{eq:intro1.1:BenOno}
   (-\Delta)^{1/2}u=-u+u^2\quad\text{ in }\R
\end{equation} 
modeling traveling wave solutions to the Benjamin-Ono equation in hydrodynamics and for the Peierls-Nabarro equation 
\begin{equation}\label{eq:intro1.1:Peierls}
   (-\Delta)^{1/2}u=\sin(\pi u)\quad\text{ in }\R
\end{equation} 
from continuum modeling of dislocations in crystals. These are two really special equations since they are ``completely integrable''. More precisely,  Amick and Toland~\cite{AmickToland1991} were able to find, with explicit expressions, all bounded solutions to 
\eqref{eq:intro1.1:BenOno}. This required $f(u)$ to be exactly $-u+u^2$ and $s=1/2$.\footnote{Their method consisted of identifying the equation as the boundary reaction problem \eqref{eq:halfhalf} which they could solve, for $f(u)=-u+u^2$, using complex analysis techniques.}
Later, in \cite{Toland1997} Toland also found all bounded solutions to the Peierls-Nabarro equation using the fact that a derivative of any solution is the difference of two solutions to the Benjamin-Ono equation. 

Both equations \eqref{eq:intro1.1:BenOno} and \eqref{eq:intro1.1:Peierls} admit infinitely many periodic solutions of different amplitude. Equation \eqref{eq:intro1.1:BenOno} admits, on the other hand, a ``ground state'' solution; that is, an even solution which is decreasing in $(0,+\infty)$. Instead, \eqref{eq:intro1.1:Peierls} admits a ``layer solution'', that is, a solution that is increasing in all of $\R$.

For nonlinearities different than the two above there are few works that classify solutions or prove symmetry and monotonicity properties. They are described in a recent paper  by Cabr\'e, Csato, and Mas~\cite{CabreCsatoMas2024}.

In \cite{CabreCsatoMas2024} we take up on the question of solutions to \eqref{eq:frac Lapl intro} being necessarily even with respect to any of their critical points. We restrict ourselves to periodic solutions to~\eqref{eq:frac Lapl intro}. The paper develops first the basic theory of periodic solutions for the fractional Laplacian. In particular, an important result in~\cite{CabreCsatoMas2024} states that the energy functional giving rise to periodic solutions to  \eqref{eq:frac Lapl intro} is different than the usual one for Dirichlet exterior boundary conditions. That is, the interactions of points (inside the domain versus outside) in the functional for periodic solutions differs from the Dirichlet interactions. The periodic functional is the one in the following Theorem~\ref{thm: frac Laplacian Intro}.
 
In the next result we establish the even symmetry and monotonicity of periodic constrained minimizers for a nonlocal Lagrangian whose first variation is the semilinear equation $(-\Delta)^s u=f(u)$.

\begin{theorem}[\cite{CabreCsatoMas2024}, Theorem 1.1]
\label{thm: frac Laplacian Intro}
Let $L>0$, $0<s<1$, $F\in C^1(\R)$, $\widetilde{F}\in C^1(\R)$, $c\in\R$, and set $f:=F'$, $\widetilde{f}:=\widetilde{F}'$. Assume that $u\in L^\infty(\R)$ is $2L$-periodic and minimizes the functional
\begin{equation}\label{energy_fac_lap}
\begin{split}
E(v):= \frac{c_s}{4}\int_{-L}^L dx\int_\R dy\,\frac{|v(x)-v(y)|^2}{|x-y|^{1+2s}}-\int_{-L}^L F(v)
\end{split}
\end{equation}
among all $2L$-periodic functions $v\in L^\infty(\R)$ satisfying the constraint
\begin{equation}
\label{thm1_constr}
\int_{-L}^L \widetilde{F}(v)=c.
\end{equation}

Then, up to a translation in the variable $x\in\R$, $u$ is an even function in $\R$ which is nonincreasing in $(0,L)$.
\end{theorem} 
 
If the conclusions of Theorem~\ref{thm: frac Laplacian Intro} hold also for all periodic solutions to \eqref{eq:frac Lapl intro} ---and not only for constrained minimizers--- is still an open problem.

The proof of the theorem is based on a not so well-known Riesz rearrangement inequality on~$\S^1$ from 1976, found independently by Baernstein and Taylor, and by Friedberg and Luttinger. 
 
It is worth pointing out that some semilinear equations $(-\Delta)^s u= f(u)$ do admit periodic constrained minimizers
which are nonconstant. Indeed, in a forthcoming work we will show their
existence for the equation
$(-\Delta)^s u=-u+u^p$ in $\R$, where $u>0$ and $p>1$ is subcritical, whenever the period $2L$ is large
enough.
On the other hand, one can show that, under no constraint,
there are no nonconstant bounded periodic minimizers ---not even local minimizers, meaning minimizers among
small periodic perturbations. That is, being periodic, such solutions must be constrained minimizers.

To conclude, we mention that in~\cite{CabreCsatoMas2024} we also prove the analogue of Theorem~\ref{thm: frac Laplacian Intro} for more general integro-differential operators.  More precisely, we establish it for integro-differential operators with symmetric kernels $K(|x-y|)$, for $x$ and $y$ in $\R$, in three different classes. Namely, we need to assume one of the following conditions on $K:(0,+\infty)\longrightarrow (0, +\infty)$:
\begin{itemize}
\item[$(i)$] either $K$ is convex,
\item[$(ii)$] ot the function $\tau>0\mapsto K(\tau^{1/2})$ 
is completely monotonic, 
\item[$(iii)$] or $K$ is nonincreasing in $(0,L)$ and vanishes in $[L,+\infty)$.
\end{itemize}
We recall that a function $J:(0,+\infty)\to\R$ is completely monotonic if
$J$ is infinitely differentiable and satisfies
$$(-1)^k\frac{d^k}{d\tau^k}J(\tau)\geq0
\qquad\text{for all $k\geq0$ and all $\tau>0$.}
$$
We recall also that, by Bernstein's theorem, a function is completely monotonic if and only if it is the Laplace transform of a nonnegative measure.

Classes $(i)$ and $(ii)$ of kernels are more interesting than the third one since they do not refer to the period $2L$, that is in principle unknown when searching for periodic solutions. The first class contains the kernel of the fractional Laplacian $K(t)=c_s t^{-1-2s}$, but not the kernel
\begin{equation}
\label{eq:intro Del kernel}
K(t)=(t^2+a^2)^{-(n+s)/2}  \qquad\text{  for some $a> 0$ and $n\geq 2$}
\end{equation} 
appearing in the study~\cite{CabreCsatoMasDelaunay} of nonlocal Delaunay cylinders. Thus, the first class does not allow to conclude symmetry for the nonlocal Delaunay cylinders. Instead, our second class includes both \eqref{eq:intro Del kernel} and also the kernel of the fractional Laplacian.

In \cite{CabreCsatoMas2024} it is shown that none of the two classes contains the other. Both, however, are included in the class of nonincreasing kernels, for which it would be tempting to conjecture that symmetry (as in Theorem~\ref{thm: frac Laplacian Intro}) also holds. However, in  \cite{CabreCsatoMas2024} we show that the theorem is false
in general for nonincreasing kernels. It remains open to find a ``natural'' class containing both $(i)$ and $(ii)$ and for which symmetry holds.

\bibliographystyle{siam}
\bibliography{cabre}

\end{document}